\documentclass{gtart}

\def\ifplaintex{\expandafter\ifx\csname documentclass\endcsname\relax}

\def\gtp{{\mathsurround=0pt\it $\cal G\mskip-2mu$eometry \&\ 
$\cal T\!\!$opology $\cal P\!$ublications}}  

\def\recd{{\small Received:\qua\receiveddate\ifx\reviseddate\relax
\else\qquad Revised:\qua\reviseddate\fi\par}} 

\def\lognumber#1{\def\thelognumber{#1}}
\def\volumenumber#1{\def\thevolumenumber{#1}}
\def\volumeyear#1{\def\thevolumeyear{#1}}
\def\papernumber#1{\def\thepapernumber{#1}}
\def\pagenumbers#1#2{\def\startpage{#1}\def\finishpage{#2}}
\def\published#1{\def\publishdate{#1}}

\def\received#1{\def\receiveddate{#1}}
\def\revised#1{\def\reviseddate{#1}}
\def\accepted#1{\def\accepteddate{#1}}

\def\asciiaddress#1{\def\theasciiaddress{#1}}

\let\\\par\let\thelognumber\relax\let\thevolumenumber\relax
\let\thepapernumber\relax\let\thevolumeyear\relax\let\startpage\relax
\let\finishpage\relax\let\publishdate\relax\let\receiveddate\relax
\let\reviseddate\relax\let\accepteddate\relax\let\theasciititle\relax
\let\theasciiauthors\relax\let\theasciiaddress\relax
\let\theasciiabstract\relax

\let\theasciiemail\relax

\ifplaintex
\font\logobig=cmssbx10 scaled 3836
\font\logomed=cmssbx10 scaled 2557
\else
\font\logobig=cmssbx10 scaled 4200
\font\logomed=cmssbx10 scaled 2800
\fi

\long\def\makeagttitle{   
\count0=\startpage
\agt\hfill      
\hbox to 45truept{\vbox to 0pt{\vglue -13truept{\logomed A\kern -.37em{\logobig 
T}\kern -.38em G}\vss}\hss}
\break
{\small Volume \thevolumenumber\ (\thevolumeyear)
\startpage--\finishpage\nl
Published: \publishdate}

\vglue .25truein

{\parskip=0pt\leftskip 0pt plus
1fil\def\\{\par\smallskip}{\Large\bf\thetitle}\par\medskip} \vglue
0.05truein

{\parskip=0pt\leftskip 0pt plus 1fil\def\\{\par}{\sc\theauthors}
\par\medskip}%
 
\vglue 0.03truein 

{\small\leftskip 25truept\rightskip 25truept{\bf Abstract}\stdspace\theabstract

{\bf AMS Classification}\stdspace\theprimaryclass
\ifx\thesecondaryclass\relax\else; \thesecondaryclass\fi\par
{\bf Keywords}\stdspace \thekeywords\par}\vglue 7truept

}   

\ifplaintex
\font\phead=cmsl9 scaled 950
\font\pnum=cmbx10 scaled 913
\font\pfoot=cmsl9 scaled 950
\headline{\vbox to 0pt{\vskip -4.5mm\line{\small\phead\ifnum
\count0=\startpage ISSN 1472-2739 (on-line) 1472-2747 (printed)
\hfill {\pnum\folio}\else\ifodd\count0\def\\{ }%
\ifx\theshorttitle\relax\thetitle\else\theshorttitle\fi\hfill{\pnum\folio}
\else\def\\{ and }{\pnum\folio}\hfill\ifx\theshortauthors\relax\theauthors
\else\theshortauthors\fi\fi\fi}\vss}}
\footline{\vbox to 0pt{\vglue 0mm\line{\small\pfoot\ifnum\count0=\startpage
\copyright\ \gtp\hfill\else
\agt, Volume \thevolumenumber\ (\thevolumeyear)\hfill\fi}\vss}}
\else
\font=cmsl9 scaled 1050
\font\lnum=cmbx10 
\font=cmsl9 scaled 1050
\makeatletter
\def\@oddhead{{\small\lhead\ifnum\count0=\startpage ISSN 1472-2739 
(on-line) 1472-2747 (printed)\hfill {\lnum\number\count0}\else\ifodd\count0
\def\\{ }\ifx\theshorttitle\relax \thetitle \else\theshorttitle\fi\hfill
{\lnum\number\count0}\else\def\\{ and }{\lnum\number\count0}
\hfill\ifx\theshortauthors\relax 
\theauthors\else\theshortauthors\fi\fi\fi}}\def\@evenhead{\@oddhead}
\def\@oddfoot{\small\lfoot\ifnum\count0=\startpage\copyright\ \gtp\hfill\else
\agt, Volume \thevolumenumber\ (\thevolumeyear)\hfill\fi}
\def\@evenfoot{\@oddfoot}
\makeatother
\fi
\let\maketitlepage\makeagttitle
\let\makeshorttitle\maketitlepage
\let\maketitle\maketitlepage

\newwrite\gtoutfile
\long\gdef\makeheadfile{  
{\def\\{, }\def\s{ }
\immediate\openout\gtoutfile head.xxx
\immediate\write\gtoutfile{To: math@arxiv.org}
\immediate\write\gtoutfile{Subject: put OR rep NNNNN:ppppp}
\immediate\write\gtoutfile{--text follows this line--}
\immediate\write\gtoutfile{Proxy-for: \ifx\theasciiauthors\relax
\theauthors\else\theasciiauthors\fi\s<\ifx\theasciiemail\relax\theemail\else\theasciiemail\fi>}
\immediate\write\gtoutfile{\noexpand\\}
\immediate\write\gtoutfile{Authors: \ifx\theasciiauthors\relax
\theauthors\else\theasciiauthors\fi}
{\def\\{ }\immediate\write\gtoutfile{Title: \ifx\theasciititle\relax
\thetitle\else\theasciititle\fi}}
\immediate\write\gtoutfile{Subj-class: GT or SG, GR etc}
\immediate\write\gtoutfile{MSC-class: \theprimaryclass\ifx\thesecondaryclass\relax\else, \thesecondaryclass\fi}
\immediate\write\gtoutfile{Journal-ref: Algebr. Geom. Topol. \thevolumenumber\s
(\thevolumeyear) \startpage-\finishpage}
\immediate\write\gtoutfile{Comments: Published by Algebraic and
Geometric Topology at}
\immediate\write\gtoutfile{\s\s\s  http://www.maths.warwick.ac.uk/agt/AGTVol\thevolumenumber/agt-\thevolumenumber-\thepapernumber.abs.html}
\immediate\write\gtoutfile{\noexpand\\}
\immediate\write\gtoutfile{}
\ifx\theasciiabstract\relax
\immediate\write\gtoutfile{\theabstract}\else
\immediate\write\gtoutfile{\theasciiabstract}\fi
\immediate\write\gtoutfile{}
\immediate\write\gtoutfile{\noexpand\\}
\immediate\write\gtoutfile{}
\immediate\closeout\gtoutfile}}  

\def\maketitlepage{\makeagttitle\makeheadfile}
\let\makeshorttitle\maketitlepage
\let\maketitle\maketitlepage

\def\ifplaintex{\expandafter\ifx\csname documentclass\endcsname\relax}

\def\gtp{{\mathsurround=0pt\it $\cal G\mskip-2mu$eometry \&\ 
$\cal T\!\!$opology $\cal P\!$ublications}}  

\def\recd{{\small Received:\qua\receiveddate\ifx\reviseddate\relax
\else\qquad Revised:\qua\reviseddate\fi\par}} 

\def\lognumber#1{\def\thelognumber{#1}}
\def\volumenumber#1{\def\thevolumenumber{#1}}
\def\volumeyear#1{\def\thevolumeyear{#1}}
\def\papernumber#1{\def\thepapernumber{#1}}
\def\pagenumbers#1#2{\def\startpage{#1}\def\finishpage{#2}}
\def\published#1{\def\publishdate{#1}}

\def\received#1{\def\receiveddate{#1}}
\def\revised#1{\def\reviseddate{#1}}
\def\accepted#1{\def\accepteddate{#1}}

\def\asciiaddress#1{\def\theasciiaddress{#1}}

\let\\\par\let\thelognumber\relax\let\thevolumenumber\relax
\let\thepapernumber\relax\let\thevolumeyear\relax\let\startpage\relax
\let\finishpage\relax\let\publishdate\relax\let\receiveddate\relax
\let\reviseddate\relax\let\accepteddate\relax\let\theasciititle\relax
\let\theasciiauthors\relax\let\theasciiaddress\relax
\let\theasciiabstract\relax

\let\theasciiemail\relax

\ifplaintex
\font\logobig=cmssbx10 scaled 3836
\font\logomed=cmssbx10 scaled 2557
\else
\font\logobig=cmssbx10 scaled 4200
\font\logomed=cmssbx10 scaled 2800
\fi

\long\def\makeagttitle{   
\count0=\startpage
\agt\hfill      
\hbox to 45truept{\vbox to 0pt{\vglue -13truept{\logomed A\kern -.37em{\logobig 
T}\kern -.38em G}\vss}\hss}
\break
{\small Volume \thevolumenumber\ (\thevolumeyear)
\startpage--\finishpage\nl
Published: \publishdate}

\vglue .25truein

{\parskip=0pt\leftskip 0pt plus
1fil\def\\{\par\smallskip}{\Large\bf\thetitle}\par\medskip} \vglue
0.05truein

{\parskip=0pt\leftskip 0pt plus 1fil\def\\{\par}{\sc\theauthors}
\par\medskip}%
 
\vglue 0.03truein 

{\small\leftskip 25truept\rightskip 25truept{\bf Abstract}\stdspace\theabstract

{\bf AMS Classification}\stdspace\theprimaryclass
\ifx\thesecondaryclass\relax\else; \thesecondaryclass\fi\par
{\bf Keywords}\stdspace \thekeywords\par}\vglue 7truept

}   

\ifplaintex
\font\phead=cmsl9 scaled 950
\font\pnum=cmbx10 scaled 913
\font\pfoot=cmsl9 scaled 950
\headline{\vbox to 0pt{\vskip -4.5mm\line{\small\phead\ifnum
\count0=\startpage ISSN 1472-2739 (on-line) 1472-2747 (printed)
\hfill {\pnum\folio}\else\ifodd\count0\def\\{ }%
\ifx\theshorttitle\relax\thetitle\else\theshorttitle\fi\hfill{\pnum\folio}
\else\def\\{ and }{\pnum\folio}\hfill\ifx\theshortauthors\relax\theauthors
\else\theshortauthors\fi\fi\fi}\vss}}
\footline{\vbox to 0pt{\vglue 0mm\line{\small\pfoot\ifnum\count0=\startpage
\copyright\ \gtp\hfill\else
\agt, Volume \thevolumenumber\ (\thevolumeyear)\hfill\fi}\vss}}
\else
\font=cmsl9 scaled 1050
\font\lnum=cmbx10 
\font=cmsl9 scaled 1050
\makeatletter
\def\@oddhead{{\small\lhead\ifnum\count0=\startpage ISSN 1472-2739 
(on-line) 1472-2747 (printed)\hfill {\lnum\number\count0}\else\ifodd\count0
\def\\{ }\ifx\theshorttitle\relax \thetitle \else\theshorttitle\fi\hfill
{\lnum\number\count0}\else\def\\{ and }{\lnum\number\count0}
\hfill\ifx\theshortauthors\relax 
\theauthors\else\theshortauthors\fi\fi\fi}}\def\@evenhead{\@oddhead}
\def\@oddfoot{\small\lfoot\ifnum\count0=\startpage\copyright\ \gtp\hfill\else
\agt, Volume \thevolumenumber\ (\thevolumeyear)\hfill\fi}
\def\@evenfoot{\@oddfoot}
\makeatother
\fi
\let\maketitlepage\makeagttitle
\let\makeshorttitle\maketitlepage
\let\maketitle\maketitlepage

\newwrite\gtoutfile
\long\gdef\makeheadfile{  
{\def\\{, }\def\s{ }
\immediate\openout\gtoutfile head.xxx
\immediate\write\gtoutfile{To: math@arxiv.org}
\immediate\write\gtoutfile{Subject: put OR rep NNNNN:ppppp}
\immediate\write\gtoutfile{--text follows this line--}
\immediate\write\gtoutfile{Proxy-for: \ifx\theasciiauthors\relax
\theauthors\else\theasciiauthors\fi\s<\ifx\theasciiemail\relax\theemail\else\theasciiemail\fi>}
\immediate\write\gtoutfile{\noexpand\\}
\immediate\write\gtoutfile{Authors: \ifx\theasciiauthors\relax
\theauthors\else\theasciiauthors\fi}
{\def\\{ }\immediate\write\gtoutfile{Title: \ifx\theasciititle\relax
\thetitle\else\theasciititle\fi}}
\immediate\write\gtoutfile{Subj-class: GT or SG, GR etc}
\immediate\write\gtoutfile{MSC-class: \theprimaryclass\ifx\thesecondaryclass\relax\else, \thesecondaryclass\fi}
\immediate\write\gtoutfile{Journal-ref: Algebr. Geom. Topol. \thevolumenumber\s
(\thevolumeyear) \startpage-\finishpage}
\immediate\write\gtoutfile{Comments: Published by Algebraic and
Geometric Topology at}
\immediate\write\gtoutfile{\s\s\s  http://www.maths.warwick.ac.uk/agt/AGTVol\thevolumenumber/agt-\thevolumenumber-\thepapernumber.abs.html}
\immediate\write\gtoutfile{\noexpand\\}
\immediate\write\gtoutfile{}
\ifx\theasciiabstract\relax
\immediate\write\gtoutfile{\theabstract}\else
\immediate\write\gtoutfile{\theasciiabstract}\fi
\immediate\write\gtoutfile{}
\immediate\write\gtoutfile{\noexpand\\}
\immediate\write\gtoutfile{}
\immediate\closeout\gtoutfile}}  

\def\maketitlepage{\makeagttitle\makeheadfile}
\let\makeshorttitle\maketitlepage
\let\maketitle\maketitlepage

\lognumber{33}
\volumenumber{1}
\volumeyear{2001}
\papernumber{33}
\pagenumbers{687}{697}
\received{10 July 2001}
\revised{14 November 2001}
\accepted{16 November 2001}
\published{19 November 2001}

\usepackage{amssymb,amsmath,epsfig,pb-diagram,lamsarrow,pb-lams}

\theoremstyle{plain}

\newtheorem{theorem}{Theorem}
\newtheorem{proposition}{Proposition}[section]
\newtheorem{lemma}[proposition]{Lemma}

\theoremstyle{definition}

\theoremstyle{remark}

\newtheorem{remark}[proposition]{Remark}

\newcommand{\psdraw}[2]
         {\begin{array}{c} \hspace{-1.3mm}
	\raisebox{-4pt}{\epsfig{figure=#1.eps,width=#2}}
	\hspace{-1.9mm}\end{array}}

\def\lbl#1{\label{#1}}

\def\BZ{\mathbb Z}

\def\BR{\mathbb R}

\def\D{\Delta}
\def\K{\mathcal K}

\def\M{\mathcal M}

\def\ihs{integral homology 3-sphere}

\def\s{$ \spadesuit $}

\def\mat#1#2#3#4{\left(
\begin{matrix}
 #1 & #2  \\
 #3 & #4   
\end{matrix}
\right)}
\def\i{^{-1}}

\def\bd{\partial}

\def\th{\theta}

\def\s{\sigma}

\def\sub{\subset}

\def\ti{\widetilde}

\def\ygraph{$\mathrm{Y}$-graph}

\def\zt{\BZ [t,t\i ]}

\renewcommand{\cl}{$\clubsuit$}

\def\th{\theta}

\def\K{\mathcal K}

\def\L{{\mathcal L}}

\def\M{\mathcal M}

\def\bd{\partial}
\def\i{^{-1}}

\begin{document}

\title{Concordance and 1-loop clovers}

\author{Stavros Garoufalidis\\Jerome Levine}
\address{Department of Mathematics, University of Warwick\\Coventry, 
CV4 7AL, UK}
\secondaddress{Department of Mathematics, Brandeis University\\Waltham, 
MA 02254-9110, USA}

\asciiaddress{Department of Mathematics, University of Warwick\\Coventry, 
CV4 7AL, UK\\and\\Department of Mathematics,
Brandeis University\\Waltham, MA 02254-9110, USA}

\email{stavros@maths.warwick.ac.uk, levine@brandeis.edu}

\url{http://www.math.gatech.edu/\char'176stavros,
http://www.math.brandeis.edu/Faculty/jlevine/}

\primaryclass{57N10}
\secondaryclass{57M25}
\keywords{Concordance, $S$-equivalence, clovers,
finite type invariants}
\begin{abstract}
We show that surgery on a connected clover (or clasper) with at least one 
loop
preserves the concordance class of a knot. Surgery on a slightly more 
special
class of clovers preserves invertible concordance. We also show that the
converse is false. Similar results hold for clovers with at least two loops 
vs. $S$-equivalence.
\end{abstract}

\maketitle


\section{Introduction}
\lbl{sec.intro}

\subsection{History}
\lbl{sub.history}

M. Goussarov and  K. Habiro have independently studied links and 
3-manifolds from the point of view of surgery on objects
called {\em Y-graphs, claspers or clovers}, respectively by \cite{Gu,H} 
and \cite{GGP}. Following the notation of \cite{GGP},
given a pair $(M,K)$ consisting of a knot $K$ in an \ihs\ $M$,
and a clover $G \subset M-K$, surgery on the framed link associated to $G$ 
produces a new pair $(M,K)_G$. Thus, by specifying a 
class of clovers $\mathfrak{c}$ we can define an equivalence relation
(also denoted by $\mathfrak{c}$) on the set $\K\M$ of knots in \ihs s and
sometimes
on its subset $\K$ of knots in $S^3$.

It is often the case that for certain classes of clovers $\mathfrak{c}$, the 
equivalence relation is related to some natural topological 
equivalence relation. In this paper we will be particularly interested in 
{\em concordance} (in the smooth category) but will also
discuss $S$-equivalence.

We begin by discussing some known facts. Using the terminology of 
\cite{GGP}, let $\mathfrak{c}^{\D\D}$ denote the class of clovers $G \subset 
S^3 -K$ of degree $1$ (that is, the class of \ygraph s) whose leaves form a 
$0$-framed unlink which bounds disks disjoint from $G$ that intersect $K$ 
geometrically twice and algebraically zero times. Surgery on such clovers 
was called a {\em double $\D$elta} move by Naik-Stanford, who showed that

\begin{theorem}
{\rm\cite{NS}}\qua 
\lbl{thm.NS}
$\mathfrak{c}^{\D\D}$ coincides with $S$-equivalence on $\K$.
\end{theorem}

Relaxing the above condition, let $\mathfrak{c}^{\mathrm{loop}}$ denote the
class of
clovers 
$G \subset M-K$ whose leaves have zero linking number with $K$. Surgery on
such clovers was called a {\em loop move} by G.-Rozansky who showed that

\begin{theorem}
{\rm\cite{GR}}\qua
\lbl{thm.GR} 
$\mathfrak{c}^{\mathrm{loop}}$ coincides with $S$-equivalence on $\K\M$.
\end{theorem}

Let us make the following definition. If $G$ is a clover in $M-K$ and $\L$ 
a set of leaves of $G$, we say $\L$ is  {\em simple} if the elements of $\L$ bound
disks in $M$ each of which intersects $K$ exactly once but whose interiors
otherwise are disjoint from $K$, $G$ and each other.
Consider now for every non-negative integer $n$, the class $\mathfrak{c}^n$ 
of clovers $G \subset S^3-K$ whose entire set of leaves is simple, and such 
that each connected component of $G$ is a graph with at least $n$ loops 
(i.e., whose first betti number is at least $n$). Kricker and Murakami-Ohtsuki
showed that

\begin{theorem}
{\rm\cite{Kr,MO}}\qua
\lbl{thm.Kr}
$\mathfrak{c}^{2}$ implies $S$-equivalence on $\K$.
\end{theorem}

In fact, if we let $\mathfrak{c}^{\rm iv}$ denote the class of clovers 
$G$ such that each component of $G$ has
at least one internal trivalent vertex, and $G$ has a simple set of leaves 
containing one leaf from each component, then it is not hard to check that
$\mathfrak{c}^2\sub\mathfrak{c}^{\rm iv}$ and \cite{Kr,MO} actually 
proved that $\mathfrak{c}^{\rm
iv}$ implies $S$-equivalence. Combining this with a recent result of
Conant-Teichner \cite{CT} we actually have:
\begin{theorem} 
{\rm\cite{CT}}\qua
$\mathfrak{c}^{\rm iv}$ coincides with $S$-equivalence on $\K$.
\end{theorem}

\subsection{Statement of the results}
\lbl{sub.statement}

In the present paper we will prove the following results.

\begin{theorem}
\lbl{thm.GL}
$\mathfrak{c}^{1}$ implies concordance on $\K$.
\end{theorem}

An different proof of Theorem \ref{thm.GL} has been obtained 
by Conant-Teichner \cite{CT} relying on the notion of {\em grope 
cobordism}.
This result was also announced by the first author in \cite{Le2}, where an
analogous statement was proved, and our proof will follow the lines of
that argument. The result was also known to Habiro, according to private
communication.

A slight refinement of the class $\mathfrak{c}^1$ relates to a classical
refinement of concordance known as {\em invertible concordance}. Recall that 
a knot in $S^3$ is called {\em double-slice} if it can be exhibited
as the intersection of a $3$-dimensional hyperplane in $\BR^4$ with an {\em
unknotted} imbedding of $S^2$ in $\BR^4$; see e.g. \cite{Su}. Such knots 
are obviously slice, and it is shown in \cite{Su} that, for any knot $K$, the
connected sum $K \sharp (-K)$ is double-slice, where $-K$ denotes the
mirror image of $K$. On the other hand the Stevedore knot is slice but not 
double-slice (see \cite{Su}). More generally, following \cite{Su}, we say that 
$K$ is {\em invertibly concordant} to $K'$ if there is a concordance $V$ from 
$K$ to $K'$ and a concordance $W$ from $K'$ to $K$ so that if we stack $W$ on 
top of $V$, the resulting concordance from $K$ to itself is diffeomorphic to 
the product concordance $(I\times S^3 ,I\times K)$. If we write $K\le K'$, 
then $\le$ is transitive and reflexive and perhaps even a partial ordering. 
It is easy to see that  $0\le K$, where $0$ denotes the trivial knot, if and 
only if $K$ is double-slice.

Let $\mathfrak{c}^{1,\mathrm{nf}}$ denote the subclass of 
$\mathfrak{c}^1$ consisting of
clovers with
{\em no forks}--- a fork is a trivalent vertex two of whose incident edges
contain a univalent vertex. Then, we will prove:

\begin{theorem}
\lbl{thm.double}
If $G$ is a clover in the class $\mathfrak{c}^{1,\mathrm{nf}}$ and $K'$ is
obtained from $K$ by surgery on $G$ then $K\le K'$.
\end{theorem}

It is natural to ask whether the converses to Theorems \ref{thm.Kr},
\ref{thm.GL} and \ref{thm.double} are true. If that were the case, one 
could
extract from the rational functions invariants of \cite{GK} many 
concordance 
invariants of knots. It was a bit of a surprise for us to show  that the 
converses are all false. 

First of all, it will follow easily from a recent result of Livingston
that:
\begin{proposition}
\lbl{thm.5}
There are $S$-equivalent knots which are not $\mathfrak{c}^2$-equivalent.
\end{proposition}

Then we will generalize some techniques of Kricker to prove:
\begin{theorem}
\lbl{thm.6}
There are double-slice knots which are not $\mathfrak{c}^1$-equivalent to 
the
unknot.
\end{theorem} 

\begin{remark}
\lbl{rem.examples}
The proofs of Proposition \ref{thm.5} and Theorem \ref{thm.6} allow one
to easily construct specific knots with the desired properties. See 
\cite[Theorem 10.1]{Li} for knots that satisfy Proposition \ref{thm.5}.
For the $(5,2)$-torus knot $T_{5,2}$, we have that $T_{5,2} \sharp 
(-T_{5,2})$
is a knot that satisfies Theorem \ref{thm.6}.
\end{remark}

\subsection{Plan of the proof}
\lbl{sub.plan}
Theorems \ref{thm.GL} and  \ref{thm.double} follow from
an analysis of the surgery link corresponding to a clover. 

Proposition \ref{thm.5} follows easily from the fact (proven recently by
Livingston \cite{Li}, using Casson-Gordon invariants) that $S$-equivalence 
does
not imply concordance.

Theorem \ref{thm.6} follows
from the fact that under surgery on 
$\mathfrak{c}^1$-clovers, the Alexander polynomial changes under a more
restrictive
way than under a concordance.

\section{Proofs}
\lbl{sec.proofs}

\subsection{Proof of Theorem \ref{thm.GL}}
\lbl{sub.thm.GL}

Suppose that $G$ is a connected clover of class $\mathfrak{c}^1$ and $L$ is 
its
associated framed link, \cite{Gu,H,GGP}. We want to show that the knot $K'$ 
obtained from $K$ by surgery on $L$ is concordant to $K$. Note that the 
manifold $M$ obtained from $S^3$ by surgery on $L$ is diffeomorphic to 
$S^3$,
see \cite{Gu,H,GGP}.

\begin{lemma}
\lbl{lem.split}
We can express $L$ as a union of two sublinks $L'$ and $L''$ such that:
\begin{itemize}
\item $L'$ is a trivial $0$-framed link in $S^3 -K$,
\item $L''$ is a trivial $0$-framed link in $S^3$.
\end{itemize}
\end{lemma}
Assuming this lemma we can complete the proof of Theorem \ref{thm.GL} as
follows.

Consider $I\times  K\sub I\times S^3$ and $\frac{1}{2}\times 
L\sub\frac{1}{2}\times (S^3 -K)$. Consider a union of disjoint disks 
$D'$ in $\frac{1}{2}\times (S^3 -K)$ bounded by $L'$ and push their 
interiors into $[0,\frac{1}{2})\times (S^3-K)$. Also consider a union of 
disjoint disks $D''$ in $\frac{1}{2}\times S^3$ bounded by $L''$ and push 
their interiors into $(\frac{1}{2},1]\times S^3$. Now
let $X\sub I\times S^3$ be obtained from $[0,\frac{1}{2}]\times S^3$ by
removing a tubular neighborhood of $D'$ and adjoining a tubular 
neighborhood 
of $D''$. The boundary of $X$ consists of $0\times S^3$ and a copy of $M$, 
which is diffeomorphic to $S^3$. Thus $X$ is diffeomorphic to $I\times 
S^3$ ( indeed, add a $D^4$ to $X$ along $0 \times S^3$ and observe that
any two imbeddings of a $4$-disk in a fixed $4$-disk are isotopic).
Moreover $X$ contains $[0,\frac{1}{2}]\times K$, which is a concordance 
from 
$0\times K\sub 0\times S^3$ to $\frac{1}{2}\times K\sub M$, which is just 
$K'$.
\qed

\begin{proof}[Proof of Lemma \ref{lem.split}]
This is a generalization of the argument used to prove Theorem 2 in 
\cite{Le2}.
Recall (eg. from \cite[Section 2.3]{GGP}) that surgery on a clover $G$ with
$n$ edges corresponds to surgery on a link $L$ of $2n$
components. Given an orientation of the edges of $G$, we can split $L$
into the disjoint union of $n$-component sublinks $L'$ and $L''$, where 
$L'$
(resp. $L''$)
consists of the sublink
of $L$ assigned to the tails of the edges of $G$ (resp. of the heads of
the edges of $G$, together with the leaves of $G$).  
As long as we avoid assigning all three of
the components at a trivalent vertex to $L'$ or $L''$, we will have the 
desired decomposition of $L$. The corresponding conditions imposed on the
orientation of the edges of $G$ are:
\begin{enumerate}
\item No trivalent vertex is a source or a sink,
\item Every edge with a univalent vertex is oriented toward the univalent
vertex.
\end{enumerate}
These are the same conditions as (i) and (ii) in the proof of Theorem 2 in
\cite{Le2} except that we now require no trivalent sinks also. But this will
follow by the same argument as in \cite{Le2} except that we need to choose the 
orientations of the cut edges more carefully. In particular we need to avoid 
choosing the orientation of two cut edges which share a trivalent vertex so 
that they both point into that vertex. But it is not hard to see that this 
can be done.
\end{proof}

The next two remarks are an addendum to Theorem \ref{thm.GL}.
 
\begin{remark}
Observe that the sublinks $L'$ and $L''$ of $L$ which are constructed from $G$
have the same number of components, and that the linking matrix of $L$
is hyperbolic.
Lemma \ref{lem.split} is analogous to the case of a knot  
which bounds a Seifert surface with a metabolic Seifert surface. In that case,
the knot is algebraically slice, and if a metabolizer can be chosen to be bands
of the Seifert surface that form a slice link, then the knot is slice.
\end{remark}

\begin{remark}
Suppose that a knot $K'$ is obtained from the unknot $K$ by surgery on a 
connected clover of class $\mathfrak{c}^1$. It follows from Theorem 
\ref{thm.GL} that $K'$ is slice. Using the calculus of clovers, one can show 
that $K'$ is actually ribbon, as observed also by Kricker and Habiro. 
\end{remark}

\subsection{Proof of Theorem \ref{thm.double}}

We need a refinement of Lemma \ref{lem.split}. Consider a connected
clover $G$ of class $\mathfrak{c}^{1,\mathrm{nf}}$ and let $L$ be its 
associated framed link.

\begin{lemma}\lbl{lem.nf}
There is a link $\bar L$ in $S^3 -K$, Kirby equivalent to $L$ in $S^3 -K$, 
so that $\bar L$ is a union of two sublinks $\bar L' ,\bar L''$, each of 
which is trivial in $S^3 -K$.
\end{lemma}

Assuming this lemma, we finish the proof following the lines of the 
argument
following Lemma \ref{lem.split}. The only difference is that we now use 
$\bar L$ instead of $L$ and that $X'=\overline{I\times S^3 -X}$, which is also
diffeomorphic to $I\times S^3$, now also contains $[\frac{1}{2} ,1]\times K$. 
Thus $M$ splits the trivial
concordance from $K$ to itself. This, by definition, means $K\le K'$.
\qed

\begin{proof}[Proof of Lemma \ref{lem.nf}]
For each univalent vertex of $G$, there is a corresponding part of $L$ 
which looks like the left part of Figure \ref{fig.nf}.
\begin{figure}[htpb]
$$
\psdraw{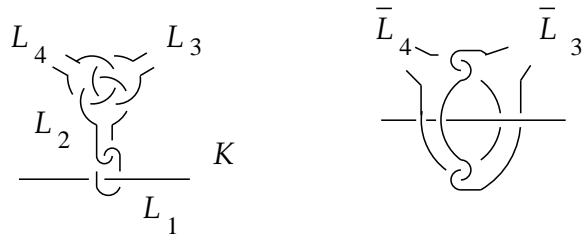}{3in}
$$
\caption{The associated link of a clover near a univalent vertex which is 
not 
a fork, before and after a Kirby move.}\lbl{fig.nf}
\end{figure}
 
Now we can perform a Kirby move (see \cite{Kr},\cite{MO}) so that the four 
component link $\{L_1,\dots,L_4\}$ in Figure \ref{fig.nf} is replaced by 
two 
component link $\{\overline{L}_3,\overline{L}_4\}$.
If we do this at every univalent vertex of $G$ we obtain the link $\bar L$.
Now consider the partition $L=L'\cup L''$ given by Lemma \ref{lem.split}.
The corresponding partition of $\overline{L}$ is given by 
$\bar L'=\{\bar K | K \in L'-\{L_1,L_2\}\}$ and
$\bar L''=\{\bar K | K \in L''-\{L_1,L_2\}\}$.
It is easy to see that both $\bar L'$ and $\bar L''$ are trivial in 
$S^3 -K$. This completes the proof.
\end{proof}

\subsection{Proof of Proposition \ref{thm.5}}
\lbl{sub.thm.5}

Assume that $S$-equivalence implies $\mathfrak{c}^2$ on $\K$. Since
$\mathfrak{c}^2$ implies
$\mathfrak{c}^1$, and $\mathfrak{c}^1$ implies concordance (by Theorem
\ref{thm.GL}),
it follows that $S$-equivalence implies concordance. This is false. 
Livingston
using Casson-Gordon invariants, shows that there are $S$-equivalent knots
which are algebraically slice, but not slice, \cite[Theorem 0.4]{Li}. 
Since Livingston uses Casson-Gordon invariants, his examples have 
nontrivial 
Alexander module.

\subsection{Proof of Theorem \ref{thm.6}}
\lbl{sub.thm.6}

We show that the Alexander polynomial $\D$ of a knot changes
in a more restrictive way under $\mathfrak{c}^1$-equivalence than under
concordance.
Recall that if $K$ and $K'$ are concordant knots, then their Alexander 
polynomials satisfy   $\D_{K'}(t)\th '(t)\th '(t\i )=\D_{K}(t)
\theta(t)\theta(t^{-1})$
for some $\theta(t),\th '(t) \in \BZ[t,t^{-1}]$ satisfying  $\theta(1)=\th
'(1)=\pm 1$. Moreover, there
are  double-slice knots with Alexander polynomial 
$\theta(t)\theta(t^{-1})$ for any such $\theta$. On the other hand,  
 
\begin{lemma} 
\lbl{lem.Kr}
Let $K$ and $K'$ be $\mathfrak{c}^1$-equivalent knots. Then,
$$
\Delta_{K'}(t)\th '(t)\th '(t\i )=\Delta_{K}(t)\theta (t)\theta (t^{-1} ) 
$$
where $\theta (t)$ and $\th '(t)$ are products of polynomials of the form 
$1\pm
t^k (t-1)^n$
for some integers $k, n$ with $n>0$.
\end{lemma}

\begin{proof}
We prove this using a generalization of an argument of Kricker \cite{Kr}.
Consider a connected clover $G$ of the class $\mathfrak{c}^1$. Suppose that
$K'$
is
obtained from $K$ by surgery on $G$. If $G$ has at least
one internal trivalent vertex, then $K$ and $K'$ are $S$-equivalent (see 
the
discussion
following Theorem \ref{thm.Kr}); in particular $\Delta_K 
(t)=\Delta_{K'}(t)$.
Otherwise, $G$ must be a {\em wheel} with a certain number $n$ of
legs and with a total of $2n$ edges. Thus, the associated
link $L'$ in $S^3 -K$ has $4n$ components (see Figure below). Using 
the Kirby move in Figure \ref{fig.nf} at every leaf of $G$ we see that 
$L'$ is Kirby-equivalent in $S^3 -K$ to a link $L$ with $2n$ components, 
whose
components can be numbered in pairs $l_1 ,r_1 ,\dots ,l_n ,r_n$ so that:
\begin{enumerate}
\item $l_i$ (resp. $r_i$) bounds a disk $d_i$ (resp. $e_i$) in $S^3 -K$,
\item $d_i\cap e_i$, for $1\le i\le n$, each consists of two oppositely
oriented
clasps,
\item  $e_i\cap d_{i+1}$, for $1\le i<n$ and $e_n\cap d_1$ each consists of 
a
single clasp, and
\item there are no other intersections among the disks.
\end{enumerate}
An example for $n=2$ is shown below:
$$
\psdraw{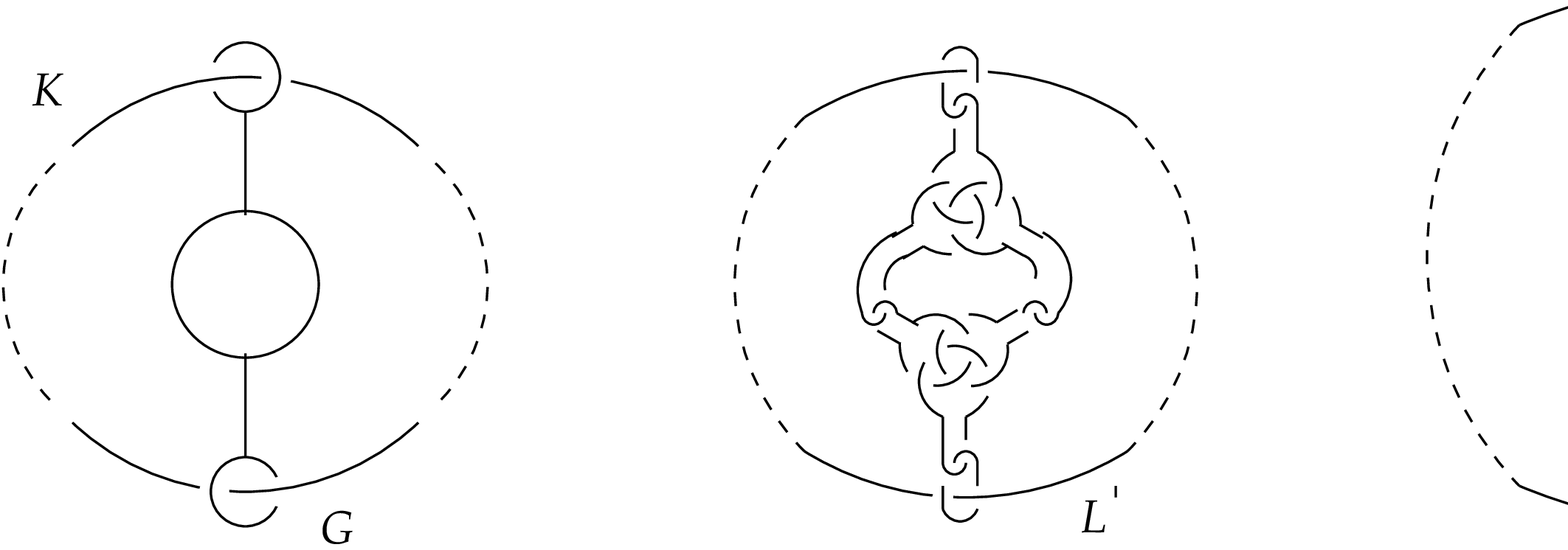}{4in} 
$$
We can now lift $d_i$ and $e_i$ to disks, $\ti d_i$ and $\ti e_i$, in the
infinite cyclic cover $\ti X$ of $X=S^3
-K$. The lifts of $l_i ,r_i$ form a link $\ti L$ in  $\ti X$ which has a
linking
matrix $B$ with entries in $\BZ [t,t\i ]$. To compute $B$ note that we can
choose the lifts $\ti d_i$ and $\ti e_i$ so that:
\begin{enumerate}
\item $\ti d_i\cap\ti e_i$ consists of a single clasp, for every $i$,
\item $\ti d_i\cap t(\ti e_i )$ consists of a single clasp, oriented 
opposite
to
that in (1), for every $i$,
\item $\ti e_i\cap\ti d_{i+1}$, for $1\le i<n$, consists of a single clasp, 
and
\item $\ti e_n\cap t^k (\ti d_1 )$, for some integer $k$,  consists of a 
single
clasp.
\end{enumerate}
In (4), $k$ (up to sign) is just the linking number of $K$ with the 
imbedded
wheel of $G$.

Now it follows from this intersection data and the fact that $L$ is 
$0$-framed
that we can orient $L$ so that the linking matrix $B$ is given by 
$$
B=\mat 0 D {D^\star} 0
\hspace{.5cm} \text{ where } \hspace{.5cm}
D=\left( \begin{array}{ccccc}t-1&1&0&\ldots&0\\
0&t-1&1&0&\vdots\\
\vdots&\ddots&\ddots&\ddots&0\\
0&\ldots&0&t-1&1\\
\pm t^k&0&\ldots&0&t-1
\end{array}\right) .
$$
For any matrix $A$ over $\BZ [t,t\i ]$, $A^\star$ denotes the conjugate
(under the involution $t\leftrightarrow t\i$) transpose of $A$.
The desired result $\D_{K'}(t)=\D_K
(t)\th (t)\th (t\i )$ is now a consequence of the following lemma, which is
proved by a standard argument going back to
Kervaire-Milnor, generalized to covering spaces (see for example
\cite[p.140]{Le1}). 
\end{proof}

Suppose $K\sub S^3$ is a knot, $L$ a framed link in $X=S^3 -K$, and $K'\sub
S^3_L$ the knot produced from $K$ by surgery on $L$. Assume that the 
components of $L$ are null-homologous in $X$ and the components of 
$\ti L\sub\ti X$, the lift
of $L$ into $\ti X$, are null-homologous. In this case we have well-defined
linking numbers of the components of $\ti L$ which are organized into a matrix
$B$ with entries in $\zt$ in the usual way. Let $A(K)=H_1(\ti X)$ and $A(K')=
H_1(\ti Y)$ denote the {\em Alexander modules} of $K, K'$, where $Y=S^3_L-K'$.

\begin{lemma}\lbl{lem.surg}
There is an exact sequence of $\zt$-modules
$$ 0\to M\to A(K')\to A(K)\to 0$$
where $M$ is a module with presentation matrix $B$. In particular,
$\D_{K'}=\D_K \det(B)$.
\end{lemma}
\begin{proof} Observe that $\ti Y=\ti X_{\ti L}$. 
Consider the following diagram of exact sequences of
$\zt$-modules.
{\small$$
\divide\dgARROWLENGTH by2
\begin{diagram}
\node[3]{H_2 (\ti Y, \ti X \mskip-1.9mu-\mskip-1.9mu \ti L )}\arrow{s,r}{\bd_*}\\
\node{H_2 (\ti X)}\arrow{e}\node{H_2(\ti X,\ti X \mskip-1.9mu-\mskip-1.9mu\ti
L)}\arrow{e}\node{H_1 (\ti X\mskip-1.9mu-\mskip-1.9mu\ti L)}\arrow{s}\arrow{e,t}{i_*}\node{H_1 (\ti
X)}\arrow{e}\node{H_1 (\ti X,\ti X\mskip-1.9mu-\mskip-1.9mu\ti L)}\\
\node[3]{H_1 (\ti Y)}\arrow{s}\\
\node[3]{H_1 (\ti Y,\ti X\mskip-1.9mu-\mskip-1.9mu\ti L)}
\end{diagram}
$$}%
Notice that $H_1 (\ti X,\ti X-\ti L)=H_1 (\ti Y,\ti X-\ti L)=0$. Moreover,
$H_2 (\ti X,\ti X-\ti L)$ is freely generated by the meridian
disks of $L$, lifted to $\ti X$, and $H_2 (\ti Y,\ti X-\ti L)$ is freely
generated by the disks attached by the surgeries. Thus, since the components of
$\ti L$ are null-homologous in $\ti X$, $i_*\circ\bd_* =0$. Also note that $H_2
(\ti X)=0$ and so we have a mapping 
$$H_2 (\ti Y,\ti X-\ti L)\to H_2 (\ti X,\ti X-\ti L)$$
induced by $\bd_*$, which can be interpreted as expressing the longitudes of
$\ti L$ as linear combinations of the meridians of $\ti L$ in $H_1 (\ti X-\ti
L)$. Therefore this map is given by the linking numbers of $\ti L$ and has $B$
as a representative matrix. This completes the proof of Lemma
\ref{lem.surg} and, as a consequence,  Lemma \ref{lem.Kr}.
\end{proof}

To complete the proof of  Theorem \ref{thm.6}  we need the
following lemma. 
\begin{lemma}\lbl{lem.root}
Let $f(t)$ be a polynomial of the form $1\pm t^k (t-1)^n$, for any integers
$k,n$ with $n\not= 0$. Then any root of $f(t)$ which lies on the unit 
circle
must be of the
form $e^{\pm\pi i/3}$. 
\end{lemma}
\begin{proof} If $z$ is a root of $f(t)$ then $|z|^k |z-1|^n =1$. Thus we 
have
$|z|=|z-1|=1$, from which the conclusion follows.
\end{proof}
Now choose some $\th (t)$ with a root on the unit circle different from
$e^{\pm\pi i/3}$ but with $\th (1)=1$---for example  any cyclotomic 
polynomial
of
composite order not equal to $6$. Let $K$ be a double-slice knot with 
Alexander polynomial $\th (t)\th (t\i )$ (see \cite[Theorem 3.3]{Su}). Then it
follows from Lemmas \ref{lem.Kr} and \ref{lem.root} that $K$ is not
$\mathfrak{c}^1$ equivalent to the trivial knot. 
\qed

We end with a remark concerning the inverse of surgery on a wheel. 
\begin{remark}
\lbl{rem.inverse}
Recall that
if a knot $K'$ is obtained from a knot $K$ by surgery on a \ygraph\ $G$, then 
there exists a \ygraph\ $G'$ such that $K$ is obtained from $K'$ by surgery 
on $G'$, see \cite[Theorem 3.2]{GGP}. Recall also that surgery on a wheel is 
described in terms of surgery on a union of \ygraph s, as explained in 
\cite[Section 2.3]{GGP}; in particular the inverse of surgery on a wheel
can be described in terms of surgery on a union of \ygraph s. One might guess 
that the inverse of surgery on a wheel can be described in terms of surgery
on a wheel. This is false, since the proof of Lemma \ref{lem.Kr} implies that
if $K'$ is obtained from $K$ by surgery on a wheel $G$, then $\D_K$ 
always divides (and it can happen that it is not equal to) $\D_{K'}$.
\end{remark}

{\bf Acknowledgements}\qua The authors were partially supported by NSF
        grants DMS-98-00703 and DMS-99-71802 respectively, and by an
        Israel-US BSF grant.

\ifx\undefined\bysame
	\newcommand{\bysame}{\leavevmode\hbox
to3em{\hrulefill}\,}
\fi

\Addresses\recd

\end{document}